\documentclass[9pt, technote, final, twocolumn]{IEEEtran}

\usepackage{enumerate}
\usepackage{amsmath, amsthm, amssymb}
\usepackage{float}
\usepackage{subcaption}
\usepackage{flexisym}
\usepackage{mathtools}
\usepackage{breqn}
\usepackage{array}
\usepackage{tikz}
\usetikzlibrary{shapes.geometric, arrows}
\usepackage{paralist}

\newtheorem{lemma}{Lemma}
\newtheorem{theorem}{Theorem}
\newtheorem{definition}{Definition}
\newtheorem{corollary}{Corollary}
\newtheorem{assumption}{Assumption}

\newcolumntype{C}[1]{>{\centering\arraybackslash$}m{#1}<{$}}
\newlength{\mycolwd}                                         
\begin{document}
\title{Dynamic Attack Detection in Cyber-Physical Systems with Side Initial State Information}
\author{Yuan Chen, Soummya Kar, and Jos\'{e} M. F. Moura
\thanks{Yuan Chen {\{(412)-268-7103\}}, Soummya Kar {\{(412)-268-8962\}}, and Jos\'{e} M.F. Moura {\{(412)-268-6341, fax: (412)-268-3890\}} are with the Department of Electrical and Computer Engineering, Carnegie Mellon University, Pittsburgh, PA 15217 {\tt\small \{yuanche1, soummyak, moura\}@andrew.cmu.edu}}
\thanks{This material is based on research sponsored by DARPA under agreement number DARPA FA8750-12-2-0291. The U.S. Government is authorized to reproduce and distribute reprints for Governmental purposes notwithstanding any copyright notation thereon.}
\thanks{The views and conclusions contained herein are those of the authors and should not be interpreted as necessarily representing the official policies or endorsements, either expressed or implied, of DARPA or the U.S. Government.}}
\maketitle

\begin{abstract}
This paper studies the impact of side initial state information on the detectability of data deception attacks against cyber-physical systems. We assume the attack detector has access to a linear function of the initial system state that cannot be altered by an attacker. First, we provide a necessary and sufficient condition for an attack to be undetectable by any dynamic attack detector under each specific side information pattern. Second, we characterize attacks that can be sustained for arbitrarily long periods without being detected. Third, we define the zero state inducing attack, the only type of attack that remains dynamically undetectable regardless of the side initial state information available to the attack detector. Finally, we design a dynamic attack detector that detects detectable attacks. 
\end{abstract}

\section{Introduction}\label{sect: introduction}
Cyber-physical systems (CPS) monitor and regulate many critical large-scale infrastructures such as the power grid and water distribution systems. Events such as the Maroochy Shire Council Sewage control incident and the Stuxnet malware attack have brought increased awareness to the issue of securing large scale systems~\cite{CardenasOld, Cardenas}. Smaller applications such as robotic platforms and the modern commercial automobile~\cite{CarAttack} are also equipped with intercommunicating sensor, computation, and actuator components for a variety of control tasks and can fall suspect to cyber attack. 
A malicious attacker can hijack the communication channels between the sensor, computation, and actuator components, modify the data values sent between components, and manipulate the system's behavior~\cite{TeixeiraModels}.

To ensure proper operation of CPS, it is necessary to design and implement security measures against attacks. 
One important aspect of security is attack detection that allows the system to take corrective actions and mitigate damaging behavior.  Static attack detectors check the consistency of the system output at a single time step~\cite{Liu, Kosut}, but are unable to detect any attacks on the actuators since they do not consider system dynamics~\cite{Pasqualetti}.  Reference~\cite{Pasqualetti} describes dynamic attack detectors that use the system dynamics, sensing topology, and the history of actuator inputs and sensor outputs to determine whether or not a data deception attack has occurred in a given time window. There are certain attacks, called stealthy or undetectable attacks, that no dynamic detector can detect. Stealthy dynamic attacks change the system output in such a way that the output of the system could arise from the system when it is not under attack~\cite{Pasqualetti}. 

There are several methods to implement attack detection. In~\cite{Mo} and~\cite{MoWorkshop}, the authors analyze dynamic attacks that go undetected by detectors of bad data (e.g., data resulting from sensor failures) for dynamical systems with process and sensor noise. References~\cite{Fawzi} and~\cite{Shoukry} provide algorithms to both detect and reconstruct the dynamic attack. The authors of~\cite{LLiu} use sparse optimization techniques to detect and identify deception attacks in electric power systems. 
Our previous work~\cite{ChenICASSP} uses geometric control techniques to analyze the limitations of detecting sparse sensor attacks. 
A different class of attack detectors, known as active attack detectors, determine the presence of a deception attack by randomly perturbing the system's input and measuring the output~\cite{MoReplay}. 
Reference~\cite{Willsky} surveys fault detection techniques in dynamic systems that are related to attack detection in CPS.
While previous work in attack detection~\cite{TeixeiraModels, Pasqualetti, Mo, MoWorkshop, TeixeiraReveal} focuses  on detectability of attacks, this note precisely clarifies how attack detector performance is sensitive to available information (specifically initial state information) and time horizons. 

We present four main contributions. First, we derive a necessary and sufficient condition for an attack to be undetectable when the detector has side initial state information given by an uncorrupted linear function of the initial system state. When the detector has initial state information, an attack is undetectable if and only if it induces a state in the intersection of the system's weakly unobservable subspace and the null space of the side information matrix.
Second, we show that an undetectable attack can be maintained if and only if the sum of the change in state produced by the attack and the zero input evolution of the state induced by the attack belong to the system's weakly unobservable subspace. An attack that is undetectable to a certain time point may become detectable at a future time as the detector obtains new sensor measurements. Undetectable attacks that can be maintained indefinitely are a greater security concern than attacks that become detectable after a finite time period. 
Third, we introduce the zero state inducing attack that is undetectable regardless of the detector's initial state information. We show that such an attack exists if and only if the intersection of the system's output-nulling reachable subspace over one time-step and its weakly unobservable subspace is nonzero. While access to initial state information improves the performance of attack detectors, it is practically important to identify the existence of attacks that are undetectable regardless of the detector's initial state information. 
Finally, we design a dynamic attack detector that uses side initial state information, has no false alarms, and only misses undetectable attacks. 

The rest of this note is organized as follows. In Section~\ref{sect: background}, we specify the system and attack model, review attack detection, introduce side information, and formally state the problem. Section~\ref{sect: mainResults} contains our main technical contributions. 
Section~\ref{sect: proofs} gives the proofs of our main results, section~\ref{sect: numExamples} provides a numerical example illustrating the performance of detectors with side information, and we conclude in Section~\ref{sect: conclusion}.

\section{Background}\label{sect: background}

\subsection{System Model}\label{sect: model}
The cyber-physical system is modeled by
\begin{equation}\label{eqn: ssModel}
	\begin{split}
		x(k+1) = Ax(k) + \overline{B}{u}(k) + Ba(k), \\
		y(k) = Cx(k) + \overline{D}{u}(k) + Da(k),
	\end{split}
\end{equation}
where: $x\in\mathbb{R}^n$ is the system state, $y\in\mathbb{R}^p$ is the system output, $k\in\mathbb{Z}$ is the time index, $u\in\mathbb{R}^m$ is the known input, and $a(k) \in \mathbb{R}^s$ is the unknown attack. 
Since the input $u(k)$ is known, its contribution to the output $y(k)$ is also known, and therefore, $u(k)$ can be ignored. Thus, for the remainder of the paper, unless otherwise stated, we consider the case of $u(k) \equiv 0$, $\forall k = 0, 1, \dots,$ without loss of generality. Accordingly, we modify the system model to be
\begin{equation}\label{eqn: ssModel2}
	\begin{split}
		x(k+1) = Ax(k) + Ba(k), \\
		y(k) = Cx(k) + Da(k).
	\end{split}
\end{equation}
The matrices $B$ and $D$ describe the capabilities of the attacker. We provide details on the attacker in Section~\ref{sect: attackBackground}. We use the notation $\Sigma = (A, B, C, D)$ to represent the system\footnote{The term ``system'' refers to the cyber-physical system and attacker collectively. The cyber-physical system gives the $A$ and $C$ matrices of $\Sigma$, while the attacker gives the $B$ and $D$ matrices of $\Sigma$.} in equation~\eqref{eqn: ssModel2}. Throughout, we make the following assumption.
\begin{assumption}\label{a: observable}
	The pair $(A, C)$ is observable.
\end{assumption}
\noindent Equation~\eqref{eqn: ssModel2} with Assumption~\ref{a: observable} is a standard model used in the cyber-physical security literature, e.g.,~\cite{Fawzi}, ~\cite{TeixeiraReveal}. 



We consider the following sequences: the output sequence (or system output trajectory) 
\begin{equation}\label{eqn: systemOutputSequence}
	Y(T) = \left[\begin{array}{cccc} y(0)^T & y(1)^T & \cdots & y(T)^T\end{array}\right]^T, 
\end{equation}
and the unknown attack sequence
\begin{equation}
	E(T) = \left[\begin{array}{cccc} a(0)^T & a(1)^T & \cdots & a(T)^T\end{array}\right]^T,
\end{equation}
with $T \geq n-1$. 
An attack occurs when $E(T) \neq 0$. The output trajectory for the deterministic system~\eqref{eqn: ssModel} is
\begin{equation}\label{eqn: outputTrajectoryLinearFunction2}
	Y(T) = \mathcal{O}_T x(0) + \mathcal{M}_T E(T),
\end{equation}
where $x(0)$ is the system's initial state, $\mathcal{O}_T$ is the extended observability matrix,
\begin{equation}\label{eqn: extendedObsMatrix}
	\mathcal{O}_T = \left[\begin{array}{c} C \\ CA \\ \vdots \\ CA^T \end{array}\right],
\end{equation}
and $\mathcal{M}_T$ is the input-output matrix,
\begin{equation}\label{eqn: mDef}
\mathcal{M}_T = \left[\begin{array}{*{5}{@{}C{\mycolwd}@{}}} D & 0  & 0& \cdots & 0 \\ 
		CB & D  &0&\cdots & 0 \\ 
		CAB & CB & D & \cdots & 0\\
		\vdots  &\vdots &\ddots &\ddots &\vdots \\ 
		CA^{t_1}B & CA^{t_2}B & \cdots & CB & D \end{array}\right],
\end{equation}
where $t_i = T - i$. In our results, we will also work with the extended controllability matrix $\mathcal{C}_T$:
\begin{equation}\label{eqn: cDef}
	\mathcal{C}_T  = \left[\begin{array}{cccc} A^TB & A^{T-1} B & \cdots & B \end{array}\right].
\end{equation}
The change in state produced by an attack $E(T)$ is $\mathcal{C}_T E(T)$. 

We now consider side initial state information. The detector knows the side initial state information
\begin{equation}\label{eqn: sideInformationStructure}
	y_\Omega = \Omega x(0),
\end{equation}
where $y_\Omega \in \mathbb{R}^q$ and $\Omega \in \mathbb{R}^{q \times n}$. We call $\Omega$ the side information matrix. The matrix $\Omega$ having full column rank corresponds to the case in which $y_\Omega$ gives full information about $x(0)$, i.e., assuming that we know $\Omega$, we can exactly determine $x(0)$ from $y_\Omega$ when $\Omega$ is full rank. The matrix $\Omega$ being the zero matrix corresponds to the case in which $y_\Omega$ gives no information about $x(0)$.

The side information $y_\Omega$ captures knowledge of the initial state $x(0)$ from the physical description of the system. For example, consider a remotely controlled vehicle whose state consists of its position and velocity. At $t=0$ the initial velocity is known to be $0$, since, by definition, the system was not running before $t=0$. We consider the initial position to be unknown since the vehicle is remotely controlled. We emphasize that the side information $y_\Omega$ does not rely on sensor measurements. For this reason, the attacker cannot modify the side information $y_\Omega$.

\subsection{Extended System Subspaces}\label{sect: extendedSubspace}
Throughout this note, we use properties of the system's extended observability and reachability subspaces (defined in~\cite{Trentelman} and~\cite{Molinari1}) to derive our results. We review their definitions here.

%

\begin{definition}[Weakly Unobservable Subspace $\mathcal{V}(\Sigma)$~\cite{Trentelman}]\label{def: weakUnobs}
	The weakly unobservable subspace of a system $\Sigma$, $\mathcal{V}(\Sigma)$, is the subspace of all $x \in \mathbb{R}^n$ such that, for a system with initial condition $x(0) = x$, there exists an input sequence $E(n-1)$ so that the output trajectory is $Y(n-1) = 0$.
\end{definition}
\noindent A state $x(0)$ belongs to the weakly unobservable subspace of $\Sigma$ if and only if there exists an input sequence $E(T)$ such that~\cite{Trentelman, Molinari1}
\[\mathcal{M}_TE(T) + \mathcal{O}_Tx(0) = 0 \text{ for any } T = 0, 1, 2, \dots\]
References~\cite{Molinari1, Trentelman, Molinari2, Rappaport} present approaches to calculate a basis for $\mathcal{V}(\Sigma)$.

Another extended system subspace of interest is the output-nulling reachable subspace over $k$ steps. 
\begin{definition}[Output-nulling Reachable Subspace $\mathcal{W}_k$~\cite{Molinari1}]\label{def: outputNull}
	The output-nulling reachable subspace over $k$ steps, $\mathcal{W}_k$, is the subspace of all states $x \in \mathbb{R}^n$ such that there exists an input (attack) sequence $E(k-1)$ that brings the system from $x(0) = 0$ to $x(k) = x$ while producing the output sequence $Y(k-1) = 0$. 
\end{definition}
\noindent The output-nulling reachable subspace over $k$ steps is the subspace of all states $x\in \mathbb{R}^n$ for which there exists $E(k-1) \in \mathbb{R}^{sk}$ such that $\mathcal{C}_{k-1} E(k-1) = x$ and $\mathcal{M}_{k-1} E(k-1) = 0$. 

\subsection{Dynamic Attack Detection: Preliminaries}\label{sect: attackBackground}
A dynamic attack detector, $\psi$, examines the system output $Y(T)$ and side initial state information $y_\Omega$ to determine whether or not an attack has occurred:
\begin{equation}\label{eqn: detectorDefSideInfo}
	\psi: \mathbb{R}^{p(T+1)} \times \mathbb{R}^q \rightarrow \left\{\text{Attack}, \text{No Attack}\right\},
\end{equation}
where $\text{``Attack''}$ means that an attack has occurred. We make the following assumptions.
\begin{assumption}
	The detector $\psi$ knows the matrices $A$ and $C$ in~\eqref{eqn: ssModel2} a priori. The detector $\psi$ does not know the matrices $B$ and $D$ in~\eqref{eqn: ssModel2} a priori. The detector $\psi$ a priori does not know $x(0)$ but knows the matrix $\Omega$ in~\eqref{eqn: sideInformationStructure}.
\end{assumption}
\noindent 
If we do not impose further restrictions on the detector, then, trivially, we can consider a detector $\psi$ that maps any input to the $\text{``Attack''}$ output. For this particular detector, every attack is detectable, but clearly this is not interesting. We restrict our focus to $\textit{consistent}$ attack detectors.
\begin{definition}[Consistent Attack Detector~\cite{Pasqualetti}]\label{def: consistency}
	An attack detector $\psi$ is consistent if $\psi\left(\mathcal{O}_T \theta, \Omega \theta\right) = \text{No Attack}$ for all $\theta \in \mathbb{R}^n$.
\end{definition}
\noindent Consistency is a desired property of attack detectors: consistent attack detectors do not produce false alarms. Another desired property of attack detectors is soundness.
\begin{definition}[Sound Attack Detector]\label{def: optimality}
	A consistent attack detector $\psi$ is sound if $\psi\left( Y(T), y_\Omega\right) = \text{No Attack}$ for some $Y(T)$ and $y_\Omega$, then, for any other consistent detector $\widetilde{\psi}$, $\widetilde{\psi}\left(Y(T), y_\Omega\right) = \text{No Attack}$. 
\end{definition}
\noindent An sound consistent detector is one that detects all possible attacks without violating the consistency property.

We now provide assumptions on the attacker. 
\begin{assumption}\label{a: bdInjective}
	The matrix $\left[\begin{array}{c} B \\ D \end{array}\right]$ is injective\footnote{If this matrix is not injective, we can remove the redundant columns to construct an injective matrix. In doing so, we do not change the capabilities of the attacker. Thus, this assumption is made without loss of generality.}.
\end{assumption}
\begin{assumption}
	The attacker knows the matrices $A, B, C, D$ and $\Omega$ and the system initial state $x(0)$ a priori.
\end{assumption}
\begin{assumption}
	The attacker cannot modify $y_\Omega$. 
\end{assumption}

Let $E(T)$ be an attack, let $Y(T)$ be the output of the system $\Sigma$ under attack $E(T)$, and let $y_\Omega$ be the side initial state information. Considering only consistent detectors, we define undetectable attacks as follows:
\begin{definition}[Undetectable Attack]\label{def: undetectableAttack}
	An attack $E(T)$ is undetectable if, for every consistent detector $\psi$ and any $x(0) \in \mathbb{R}^n$, $\psi\left(Y(T), y_\Omega\right) = \text{No Attack},$ where $Y(T) = \mathcal{O}_T x(0) + \mathcal{M}_T E(T)$. 
\end{definition}
\noindent A detectable attack is any attack that is not undetectable. 
We partition the set of all possible attacks (including $E(T) = 0$), $\mathbb{R}^{s(T+1)}$, into a set of undetectable attacks and a set of detectable attacks. 
\begin{definition}[Set of Undetectable Attacks $\mathcal{U}^{\Omega, T}$]\label{def: setOfUndetectable}
	The set $\mathcal{U}^{\Omega, T}$ is the union of set of all attacks $E(T) \in \mathcal{R}^{s(T+1)}$ such that $E(T)$ is undetectable and the set that only contains $E(T) = 0$.
\end{definition}
\noindent When the system is not under attack (i.e., $E(T) = 0$), consistent detectors report $\text{``No Attack''}$, so $0 \in \mathcal{U}^{\Omega, T}$.


Define an extension of an attack as follows:
\begin{definition}[Extension of an Attack]\label{def: extend}
An extension of $E(T)$, $E(T) \neq 0$, is an attack of the form
\begin{equation}
	\widehat{E}(T') = \left[\begin{array}{cccc} E(T)^T & a(T+1)^T & \cdots & a(T')^T \end{array}\right]^T,
\end{equation}
for $T' > T$.
\end{definition}
\noindent 
The attack sequence $a(T+1), \dots, a(T')$ is allowed to be the zero sequence. We provide a necessary and sufficient condition for which an undetectable attack $E(T)$ has undetectable extensions $\widehat{E}(T')$ for all $T' > T$ so that the attack sequence never becomes detectable (even as the attack detector obtains new sensor measurements at each time step). If $E(T)$ does not have an undetectable extension for \emph{all} times $T'> T$, then, at some time $T' > T$, regardless of the attack sequence $a(T+1), \dots, a(T')$, $\widehat{E}(T')$ is detectable. 

Reference~\cite{Pasqualetti} provides a necessary and sufficient condition for an attack sequence $E(T)$ to be undetectable when $\Omega = 0$.
\begin{lemma}[\cite{Pasqualetti}]\label{lem: undetectableCondition} The attack $E(T)$ is undetectable if and only~if \[\mathcal{O}_T x(0) + \mathcal{M}_T E(T) = \mathcal{O}_T x'(0)\] for some initial states $x(0), x'(0) \in \mathbb{R}^n$.
\end{lemma}
\noindent One particular form of attack that is undetectable against systems with no side initial state information is known as the zero dynamics attack.
\begin{definition}[Zero Dynamics Attack~\cite{TeixeiraModels}]\label{def: zeroDynamics}
	A zero dynamics attack is an attack $E(T) = \left[\begin{array}{ccc} a(0)^T & \cdots & a(T)^T\end{array}\right]^T$ with
\begin{equation}\label{eqn: zeroDynamicsIndividualAttack}
	a(k) = \lambda^k g,
\end{equation}
where $g \neq 0$ and $\lambda \in \mathbb{C}$ satisfy
\begin{equation}\label{eqn: zeroDynamicsMatrixCondition}
	\left[\begin{array}{cc} \lambda I - A & -B \\ C & D \end{array}\right]\left[\begin{array}{c} \theta \\ g \end{array}\right] = 0.
\end{equation}
\end{definition}
\noindent A zero dynamics attack exists if and only if there exists $\lambda \in \mathbb{C}$ for which there is a nonzero solution to~\eqref{eqn: zeroDynamicsMatrixCondition}~\cite{TeixeiraModels, Pasqualetti}. Since, by Assumption~\ref{a: bdInjective}, the matrix $\left[\begin{array}{cc} B^T & D^T\end{array}\right]^T$ is injective, and $g\neq 0$, we have that $\theta \neq 0$. By construction, a zero dynamics attack satisfies \[\mathcal{M}_T E(T) + \mathcal{O}_T \theta = 0.\] Therefore, a zero dynamics attack satisfies the condition given in Lemma~\ref{lem: undetectableCondition}, where $\theta = x(0) - x'(0)$. We consider $T \geq n-1$, so $\mathcal{O}_T$ is injective since $(A, C)$ is observable. Since $\theta \neq 0$, a zero dynamics attack produces a nonzero change to the output of the system. Zero dynamics attacks are also related to malicious attacks against distributed function calculation~\cite{Sundaram}. 

We introduce the zero state inducing attack:
\begin{definition}[Zero State Inducing Attack]\label{def: zeroState}
An attack sequence $E(T)$ is called a zero state inducing attack if it satisfies $\mathcal{M}_T E(T) = 0$. 
\end{definition}
\noindent The name zero-state inducing attack refers to the property that such an attack does not change the system sensor output, i.e., the change in output is equal to the response of the system when its initial state is $x(0) = 0$. We show that the zero state inducing attack is undetectable regardless of the detector's side information matrix $\Omega$. It is the only type of attack to remain undetectable even if $\Omega$ is full rank.

\subsection{Problem Statement}\label{sect: statement}
Consider a system $\Sigma = (A, B, C, D)$ over a time interval $0, 1, \dots, T$, $T \geq n-1$, with initial state $x(0)$ and side initial state information $y_\Omega = \Omega x(0)$. We consider the following four main problems:
\begin{inparaenum}[1)] 
	\item find the set of all undetectable attacks, $\mathcal{U}^{\Omega, T}$ ;
	\item determine which attacks $E(T) \in \mathcal{U}^{\Omega, T}$ have undetectable extensions up to any time $T' > T$;
	\item determine if there exists an arbitrarily long zero state inducing attack against $\Sigma$ and;
	\item design a consistent detector that uses side information and detects all detectable attacks.
\end{inparaenum}

\section{Main Results}\label{sect: mainResults}
\subsection{Initial State Information and Undetectable Attacks}
First, we find a necessary and sufficient condition for an attack to be undetectable, when the attack detector has side initial state information $y_\Omega$. Let $\mathcal{N}(\Omega)$ be the null space of $\Omega$. 

\begin{theorem}[Undetectable Attacks with Side Initial State Information]\label{thm: undetectablePartial}
	An attack $E(T)$ is undetectable ($E(T) \in \mathcal{U}^{\Omega, T}$) if and only if there exists $\theta \in \mathcal{N}\left(\Omega\right)\cap\mathcal{V}(\Sigma)$ for which $\mathcal{M}_T E(T) = -\mathcal{O}_T \theta$. 
\end{theorem}

\noindent Theorem~\ref{thm: undetectablePartial} states that an attack $E(T)$ is undetectable over the time interval $0, \dots, T$ if and only if the output contributed by the attack (i.e., $\mathcal{M}_T E(T)$) equals the negative of the output of the system operating without attack from an initial state $\theta$, where $\theta$ belongs to the intersection of the system's weakly unobservable subspace, $\mathcal{V}(\Sigma)$, and the null-space of the side information matrix, $\mathcal{N}(\Omega)$. We call $\theta$ the state induced by the attack. If $\mathcal{N}(\Omega)$ has dimension strictly less than $n$ (i.e., if the side initial state information is non-trivial), then, by using the side initial state information $y_\Omega$, an attack detector may be able to detect attacks that would otherwise be undetectable (in the absence of side information).

Theorem~\ref{thm: undetectablePartial} is valid for any side information matrix $\Omega$. 

\begin{corollary}[No Initial State Information: $\Omega = 0$]\label{cor: timeZeroAttacks}
	An attack $E(T)$ is undetectable if and only if $\mathcal{M}_T E(T)$ = $-\mathcal{O}_T \theta$ for some $\theta \in \mathcal{V}(\Sigma)$ when $\Omega = 0$.
\end{corollary}

\noindent By construction, a zero dynamics attack $E(T)$ satisfies $\mathcal{M}_TE(T) + \mathcal{O}_T \theta = 0$, where $\theta \neq 0$ and $g\neq 0$ (which is used to define $E(T)$) is a solution to equation~\eqref{eqn: zeroDynamicsMatrixCondition}. There may be other undetectable attacks aside from zero dynamics attacks when $\Omega = 0$. 

\begin{corollary}[Full Initial State Information: $\Omega$ has full column rank]\label{cor: fullInitStateAttacks}
	An attack $E(T)$ is undetectable if and only if $\mathcal{M}_TE(T) = 0$ when $\Omega$ has full column rank.
\end{corollary}
\noindent According to Corollary~\ref{cor: fullInitStateAttacks}, the only type of attack that is undetectable when the initial state is completely known to the detector is the zero state inducing attack. Figure 1 illustrates the results of Theorem~\ref{thm: undetectablePartial} and its corollaries. Undetectable attacks presented in the literature~\cite{Pasqualetti, Fawzi, Shoukry} rely on the fact that the initial state is unknown to the detector in order to be stealthy. As Theorem~\ref{thm: undetectablePartial} and Figure 1 show, however, that even when the detector knows the initial state completely, there may still be undetectable attacks. For the special case of $\Omega = 0$, Theorem~\ref{thm: undetectablePartial} is consistent with the results presented in~\cite{Pasqualetti}. 

\begin{figure}[h!]\label{fig: sideInfoResults}
	\centering
	\begin{subfigure}[h!]{\linewidth}
	\centering
	\begin{tikzpicture}[scale = .93]
		\draw[thick, rounded corners = 10pt] (-3, -1) rectangle (3, 1);
		\node [above left] at (3, -1) {$\mathbb{R}^{s(T+1)}$};
		\draw [ultra thick, rounded corners = 6 pt] (-2.6, .85) rectangle (.1, -.85);
		\node [above right] at (-2.5, -.85) {$\mathcal{U}^{\Omega, T}$};
		\node [rectangle, rounded corners = 6pt, draw = black, thin, minimum width = 1cm, minimum height = 1 cm] at (-2, .25) {$ZS$};
		\node [rectangle, rounded corners = 6pt, draw = black, thin, minimum width = 1cm, minimum height = 1 cm] at (-.5, .25) {$ZD$};
	\end{tikzpicture}
	\caption{$\Omega = 0$}
	\label{fig: sideInfoA}
	\end{subfigure}
	\vspace{\baselineskip}

	\begin{subfigure}[h!]{\linewidth}
	\centering
	\begin{tikzpicture}[scale = .93]
		\draw[thick, rounded corners = 10pt] (-3, -1) rectangle (3, 1);
		\node [above left] at (3, -1) {$\mathbb{R}^{s(T+1)}$};
		\draw [ultra thick, rounded corners = 6 pt] (-2.6, .85) rectangle (-.5, -.85);
		\node [above right] at (-2.5, -.85) {$\mathcal{U}^{\Omega, T}$};
		\node [rectangle, rounded corners = 6pt, draw = black, thin, minimum width = 1cm, minimum height = 1 cm] at (-2, .25) {$ZS$};
		\node [rectangle, rounded corners = 6pt, draw = black, thin, minimum width = 1cm, minimum height = 1 cm, fill = white] at (-.5, .25) {$ZD$};
	\end{tikzpicture}
	\caption{$\Omega \neq 0$, $\Omega$ is \textit{not} full rank}
	\label{fig: sideInfoB}
	\end{subfigure}
	\vspace{\baselineskip}

	\begin{subfigure}[h!]{\linewidth}
	\centering
	\begin{tikzpicture}[scale = .93]
		\draw[thick, rounded corners = 10pt] (-3, -1) rectangle (3, 1);
		\node [above left] at (3, -1) {$\mathbb{R}^{s(T+1)}$};
		\draw [ultra thick, rounded corners = 6 pt] (-2.5, .75) rectangle (-1.5, -.25);
		\node [above right] at (-2.5, -.85) {$\mathcal{U}^{\Omega, T}$};
		\node [rectangle, rounded corners = 6pt, draw = black, thin, minimum width = 1cm, minimum height = 1 cm] at (-2, .25) {$ZS$};
		\node [rectangle, rounded corners = 6pt, draw = black, thin, minimum width = 1cm, minimum height = 1 cm] at (-.5, .25) {$ZD$};
	\end{tikzpicture}
	\caption{$\Omega$ is full rank}
	\label{fig: sideInfoC}
	\end{subfigure}
	
	\caption{The set of all undetectable attacks $\mathcal{U}^{\Omega, T}$ depends on the side initial state information available to the attack detector. $ZS$ and $ZD$ are the set of all zero state inducing attacks and the set of all zero dynamics attacks, respectively.}
\end{figure}
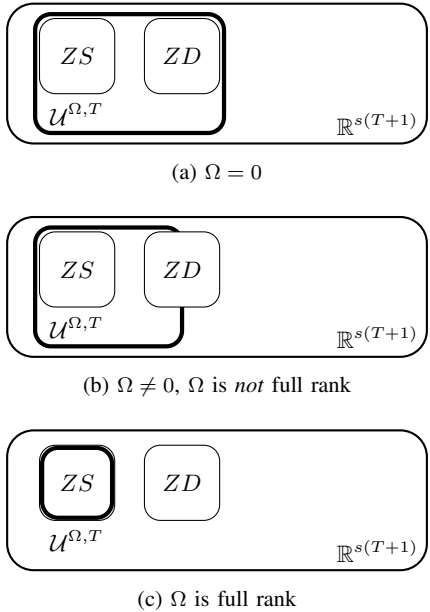

\subsection{Extensions of Undetectable Attacks}
Second, we provide a necessary and sufficient condition for an undetectable attack $E(T)$ (with $T \geq n-1$) to have an undetectable extension $\widehat{E}(T')$. Consider an attack $E(T) \in \mathcal{U}^{\Omega, T}$, $E(T) \neq~0$. 

\begin{theorem}[Extensions of Undetectable Attacks]\label{thm: undetectableExtension}
There exists an undetectable extension $\widehat{E}(T')$ of $E(T)$ for all $T' > T$ if and only if $\left(\mathcal{C}_TE(T) + A^{T+1}\theta\right) \in \mathcal{V}(\Sigma)$, where $\theta$ satisfies $\mathcal{M}_TE(T) = -\mathcal{O}_T\theta$ and $\theta \in \mathcal{N}\left(\Omega\right)\cap\mathcal{V}\left(\Sigma\right)$.
\end{theorem}

\noindent Theorem~\ref{thm: undetectableExtension} states that an undetectable attack $E(T)$ has an undetectable extension $\widehat{E}(T')$ for any $T' > T$ if and only if the sum of the change in state produced by the attack ($\mathcal{C}_T E(T)$) and the zero-input state response of the state induced by the attack ($A^{T+1} \theta$) belongs to the system's weakly unobservable subspace ($\mathcal{V}\left(\Sigma\right)$). If an attack $E(T)$ satisfies the conditions given in Theorem~\ref{thm: undetectableExtension}, then for any time $T'>T$, there exists a particular sequence of attacks $a(T+1), \dots, a(t)$ such that $\widehat{E}(T')$ is undetectable at time $T'$. Conversely, if an attack $E(T)$ does not satisfy the above condition, then at some time $T' > T$, all extensions $\widehat{E}(T')$ of $E(T)$ are detectable. In this case, all extensions $\widehat{E}(T')$ are detectable by time $T'$ because the detector obtains sensor measurements $y(T+1), \dots, y(T'+1)$ (even though $E(T)$ was undetectable). 

\subsection{Zero State Inducing Attack}\label{sect: zeroStateInducing}
Third, we provide a necessary and sufficient condition for the existence of a zero state inducing attack that can be maintained for a arbitrarily long time. 
We restrict our focus to zero state inducing attacks that begin at time $0$. 
This is to prevent trivial lengthening by appending a fixed length zero state inducing attack $E(T)$ to a zero vector\footnote{This is not a restriction on the \textit{definition} of the zero state inducing attack. An attack $E(T)$ with nonzero first attack time can still be a zero state inducing attack if $\mathcal{M}_T E(T) = 0$.}.
\begin{theorem}[Arbitrarily Long Zero State Inducing Attacks]\label{thm: arbLongZero}
	There exists an attack $E(T)$ against the system $\Sigma$ that begins at time $0$ such that $\mathcal{M}_T E(T) = 0$ for any $T = 0, 1, \dots$ if and only if $\mathcal{W}_1 \cap \mathcal{V}(\Sigma) \neq  \{0\}$, where $\mathcal{W}_1$ is the output-nulling reachable subspace over one time step.
\end{theorem}

\noindent Theorem~\ref{thm: arbLongZero} states that there exists an arbitrarily long zero state inducing attack against a system $\Sigma$ if and only if the intersection of the system's weakly unobservable subspace, $\mathcal{V}(\Sigma)$ and its output-nulling reachable subspace over one step, $\mathcal{W}_1$ is nonzero. 


\subsection{Attack Detection With Side Information}\label{sect: detectorSideInfo}
We design a consistent dynamic attack detector that detects all attacks $E(T)$ that do not belong to $\mathcal{U}^{\Omega, T}$. 
Our dynamic detector operates sequentially: at every time instant $k$ (with the exception of an initialization period), the detector collects new sensor outputs $y(k)$ and makes a decision on whether or not the system was attacked in the time period up to time $k$. 
Our detector only uses a finite window of sensor measurements in each time interval, which offers advantageous in implementation over detectors that use the entire history of sensor measurements.

First, define $\overline{Y}(k)$ as the $l$-length window of sensor measurements ending at time $k$, where $k \geq l-1$:
\begin{equation}\label{eqn: windowedY}
	\overline{Y}(k) \!\!=\!\! \left[\begin{array}{cccc} \!\!y(k-l+1)^T \!& \!y(k-l+2)^T \!&\! \cdots \!& \!\!\!y(k)^T\!\! \end{array}\right]^T.
\end{equation}
The attack detector makes a decision at every time instant starting at $l-1$. Second, define $\widehat{Y}(k)$, the input to the attack detector at time $k$, as follows:
\begin{equation}\label{eqn: hatYk}
	\widehat{Y}(k) = \left\{\begin{array}{ll} \left[\begin{array}{cc} y_\Omega^T & \overline{Y}(k)^T\end{array}\right]^T, & k = l-1 \\
				\overline{Y}(k), & k = l, l+1, \dots \end{array}\right..
\end{equation}
Third, define the orthogonal projection (operator) onto the range space of a matrix $\mathcal{K}$ (where $\mathcal{K}$ has full column rank) as
\begin{equation}\label{eqn: projector}
	\Pi_{\mathcal{K}}= \mathcal{K}\left(\mathcal{K}^T\mathcal{K}\right)^{-1} \mathcal{K}^T.
\end{equation}

We construct the detector $\psi$ as
\begin{equation}\label{eqn: windowedDetector}
	\psi\left(\widehat{Y}(k)\right) = \left\{\begin{array}{ll} \text{No Attack,} & \widehat{Y}(k) = \Pi_{\mathcal{K}(k)}\widehat{Y}(k) \\
						\text{Attack,} & \text{Otherwise}\end{array}\right.,
\end{equation}
where
\begin{equation}\label{eqn: K_k}
	\mathcal{K}(k) = \left\{\begin{array}{ll} \left[\begin{array}{cc} \Omega^T & \mathcal{O}_{l-1}^T\end{array}\right]^T, & k = l-1 \\
					\mathcal{O}_{l-1}, & k = l, l+1, \dots, \end{array}\right..
\end{equation}
The detector decides that no attack has occurred in the time interval $0, \dots, T$ if $\psi\left(\widehat{Y}(l-1)\right) = \psi\left(\widehat{Y}(l)\right) = \cdots = \psi\left(\widehat{Y}(T)\right) = \text{No Attack}$. 

\begin{theorem}[Consistency and Soundness of $\psi$]\label{thm: windowedDetectionOptimal}
	For $l \geq n+1$, where $n$ is the dimension of the system state space, $\psi\left(\widehat{Y}(l-1)\right) = \psi\left(\widehat{Y}(l)\right) = \cdots = \psi\left(\widehat{Y}(T)\right) = \text{No Attack}$ if and only if $Y(T) = \mathcal{O}_T x(0)$ and $y_\Omega = \Omega x(0)$ for some $x(0) \in \mathbb{R}^n$. 
\end{theorem}
\noindent The detector $\psi$ is consistent and sound when the window length $l$ is sufficiently long. The novelty of our detector is its use of the available side information $y_\Omega$. Detectors that do not use side information (e.g., fault detectors such as those presented in~\cite{Willsky}) may still detect some attacks, but, following Theorem~\ref{thm: undetectablePartial}, such detectors may not be sound. That is, there are certain attacks that are only detectable if the detector uses side information $y_\Omega$.

\section{Proof of Main Results}\label{sect: proofs}

\subsection{Proof of Theorem~\ref{thm: undetectablePartial}}\label{sect: undetectablePartial}
First, we provide an intermediate result by modifying Lemma~\ref{lem: undetectableCondition} to account for attack detectors with side information $y_\Omega$. Consider a system $\Sigma = (A, B, C, D)$ equipped with an attack detector that has side information matrix $\Omega$.

\begin{lemma}\label{lem: undetectableConditionSide}
	An attack $E(T)$ against the system $\Sigma$ is undetectable if and only if $\mathcal{M}_T E(T) + \mathcal{O}_T x(0) = \mathcal{O}_T x'(0)$ and  $\Omega x(0)$ = $\Omega x'(0)$ for some initial states $x(0), x'(0) \in \mathbb{R}^n$.
\end{lemma}

\noindent We use the above Lemma to prove Theorem~\ref{thm: undetectablePartial}

\begin{proof}[Proof (Theorem~\ref{thm: undetectablePartial})]
\textit{(If)} Let $x(0)$ be the initial state of the system. Let $E(T)$ be an attack such that $\mathcal{M}_T E(T) = -\mathcal{O}_T\theta$ for $\theta \in \mathcal{N}\left(\Omega\right) \cap \mathcal{V}\left(\Sigma\right)$. Let $x'(0) = x(0) - \theta$. Then $\mathcal{M}_TE(T) + \mathcal{O}_Tx(0) = \mathcal{O}_Tx'(0)$. In addition, since $\theta \in \mathcal{N}\left(\Omega\right)$, $\Omega x'(0) = \Omega \left(x(0) - \theta\right) = \Omega x(0)$. Thus, for any $x(0)$, there exists $x'(0)$ such that $\mathcal{M}_T E(T) + \mathcal{O}_T x(0) = \mathcal{O}_T x'(0)$ and $\Omega x(0) = \Omega x'(0)$, which means, by Lemma~\ref{lem: undetectableConditionSide}, $E(T)$ is an undetectable attack. Thus, $E(T) \in \mathcal{U}^{\Omega, T}$.  

\textit{(Only If)} Let $x(0)$ be the initial state of the system. Let $E(T) \in \mathcal{U}^{\Omega, T}$. Then, by Lemma~\ref{lem: undetectableConditionSide}, there exists $x'(0) \in \mathbb{R}^n$ such that $\mathcal{M}_T E(T) + \mathcal{O}_T x(0) = \mathcal{O}_T x'(0)$ and $\Omega x(0) = \Omega x'(0)$. Let $\theta = x(0) - x'(0)$. Substituting for $\theta$ we have that $\mathcal{M}_T E(T) = -\mathcal{O}_T \theta$ and $\Omega \theta = 0$. Thus, $\mathcal{M}_T E(T) = -\mathcal{O}_T \theta$ for $\theta \in \mathcal{N}\left(\Omega\right) \cap \mathcal{V}(\Sigma)$. 
\end{proof}

\subsection{Proof of Theorem~\ref{thm: undetectableExtension}}\label{sect: undetectableAttackExtensions}

\begin{proof}
\textit{(Only If)} We show that if there exists an undetectable extension $\widehat{E}(T')$ for all $T' > T$, then, necessarily, $\left(\mathcal{C}_TE(T) + A^{T+1}\theta\right) \in \mathcal{V}(\Sigma)$. Let
\[\widehat{E}(T') = \left[\begin{array}{cccc} E(T)^T & a(T+1)^T & \cdots & a(T')^T \end{array}\right]^T \]
be an undetectable extension of $E(T)$. Since $\widehat{E}(T')$ is undetectable, then, by Theorem~\ref{thm: undetectablePartial}, it must satisfy $\mathcal{M}_{T'}\widehat{E}(T') + \mathcal{O}_{T'} \theta ' = 0$ for some $\theta ' \in \mathcal{N}\left(\Omega\right)\cap\mathcal{V}\left(\Sigma\right)$. 

We first show that $\theta ' = \theta$. We partition the matrix $\mathcal{M}_{T'}$ as follows:
\begin{equation}\label{eqn: mtPart}
	\mathcal{M}_{T'} = \left[\begin{array}{cc} \mathcal{M}_T & 0 \\ \mathcal{Q}_{T'}^T & \mathcal{M}_{T'-T-1} \end{array} \right],
\end{equation}
where $\mathcal{Q}_{T'}^T = \mathcal{O}_{T'-T-1}\mathcal{C}_T$.
Substituting for the partitioned versions of $\mathcal{M}_{T'}$ and partitioning $\mathcal{O}_{T'}$, we have
\begin{equation}\label{eqn: undetectablePartition1}
	\left[\!\!\begin{array}{ccc} \mathcal{M}_T & 0 & \mathcal{O}_T \\
		\mathcal{Q}_{T'}^T & \mathcal{M}_{T'-T-1} & \mathcal{O}_{T'-T-1}A^{T+1} \end{array}\!\!\right] \!\!\left[\!\!\begin{array}{c} \widehat{E}(T') \\ \theta ' \end{array}\!\!\right] = 0.
\end{equation}
From the first block row of equation~\eqref{eqn: undetectablePartition1}, we have $\mathcal{M}_T E(T) + \mathcal{O}_T \theta ' = 0$, and, from the definition of $E(T)$, we have $\mathcal{M}_TE(T) + \mathcal{O}_T \theta = 0$. Thus, $\mathcal{O}_T \theta ' = \mathcal{O}_T \theta$. Since $T \geq n-1$ and $\Sigma$ is observable, $\mathcal{O}_T$ is injective, and $\theta ' = \theta$.

Substituting $\theta = \theta '$, the second block row of equation~\eqref{eqn: undetectablePartition1} gives
\begin{equation}\label{eqn: undetectablePartition3}
\begin{split}
	\mathcal{O}_{T'-T-1}&\left(\mathcal{C}_T E(T) + A^{T+1}\theta\right) \\
&+ \mathcal{M}_{T'-T-1} \left[\begin{array}{c} a(T+1) \\ \vdots \\ a(T') \end{array}\right] = 0.
\end{split}
\end{equation}
Since there exists an undetectable extension $\widehat{E}(T')$ of $E(T)$ for all $T' > T$, equation~\eqref{eqn: undetectablePartition3} must be satisfied for all $T' > T$. In particular, equation~\eqref{eqn: undetectablePartition3} is true for $T' = T + n$, which shows that $\left(\mathcal{C}_TE(T) + A^{T+1}\theta\right) \in \mathcal{V}(\Sigma)$.

\textit{(If)} If $\left(\mathcal{C}_TE(T) + A^{T+1}\theta\right) \in \mathcal{V}(\Sigma)$, then, for all $T' > T$, there exists an attack sequence 
$\left[\begin{array}{ccc} a(T+1)^T & \cdots & a(T')^T \end{array}\right]^T$
such that equations~\eqref{eqn: undetectablePartition3} is satisfied. For all $T' > T$, we construct $\widehat{E}(T')$ by appending $\left[\begin{array}{ccc} a(T+1)^T & \cdots & a(T')^T \end{array}\right]^T$ to $E(T)$. By definition of $E(T)$, we have 
$\mathcal{M}_TE(T) + \mathcal{O}_T\theta = 0,$ 
where $\theta \in \mathcal{N}\left(\Omega\right) \cap \mathcal{V}\left(\Sigma\right)$. Combining this fact with~equation~\eqref{eqn: undetectablePartition3}, we see that $\left[\begin{array}{c} \widehat{E}(T') \\ \theta ' \end{array}\right]$ satisfies equation~\eqref{eqn: undetectablePartition1} with $\theta ' = \theta$. Thus, we have 
\[\mathcal{M}_{T'} \widehat{E}(T') + \mathcal{O}_{T'}\theta = 0,\]
 which shows that $\widehat{E}(T')$ is an undetectable extension of $E(T)$.
\end{proof}

\subsection{Proof of Theorem~\ref{thm: arbLongZero}}\label{sect: zeroStateAttacks}

\begin{proof}
\textit{(If)} We construct a zero state inducing attack $E(T)$ that begins at time $0$ against $\Sigma$ of arbitrary length $T$  under the condition that $\mathcal{W}_1 \cap \mathcal{V}(\Sigma) \neq \{0\}$. The initial state of the system $\Sigma$, $x(0)$, does not affect its extended observability and reachability subspaces, so, without loss of generality, let the system have initial state $x(0) = 0$. If $\mathcal{W}_1 \cap \mathcal{V}(\Sigma)\neq \{0\}$, there exists an attack $a(0) \neq 0$ such that $x(1) = Ba(0)$, $y(0) = Da(0) = 0$, and $x(1) \in \mathcal{V}(\Sigma)$. Since $x(1) \in \mathcal{V}(\Sigma)$, for any $T$, there exists a sequence of attacks $\left[\begin{array}{cccc} a(1)^T & a(2)^T & \cdots &a(T)^T\end{array}\right]^T$ such that the output $\left[\begin{array}{cccc} y(1)^T & y(2)^T & \cdots & y(T)^T \end{array}\right]^T$ is $0$. Thus, for any $T$, there exists an attack $E(T) = \left[\begin{array}{cccc} a(0)^T & a(1)^T & \cdots & a(T)^T\end{array}\right]^T$ with $a(0) \neq 0$ such that $\mathcal{M}_TE(T) = 0$.

\textit{(Only If)} We show that if there exists a zero state inducing attack that begins at time $0$ for any $T$ against the system $\Sigma$, then $\mathcal{W}_1(\Sigma) \cap \mathcal{V}(\Sigma) \neq \{0\}$. Such an attack exists for any $T$, so it exists for $T = n$. Let 
\[E(n) = \left[\begin{array}{cccc} a(0)^T & a(1)^T  & \cdots & a(n)^T \end{array}\right]^T\]
be a zero state inducing attack with $a(0) \neq 0$. Since $E(n)$ induces the zero state, we have 
$\mathcal{M}_n E(n) = 0,$
which implies that $Da(0) = 0$. Since $\left[\begin{array}{c} B \\ D \end{array}\right]$ is injective and $Da(0) = 0$, we have $x(1) = Ba(0) \neq 0$ and $x(1) \in \mathcal{W}_1$. The sequence
\[\left[\begin{array}{cccc} a(1)^T & a(2)^T & \cdots & a(n)^T\end{array}\right]^T\] 
is an input sequence over $n$ steps such that a system with state $x(1) = Ba(0)$ produces zero output over the time period $1, \dots, n$. Since such an input sequence exists, $x(1) \in \mathcal{V}\left(\Sigma\right)$ and $x(1) \in \mathcal{W}_1\cap\mathcal{V}(\Sigma).$ 
Since $x(1) \neq 0$, $\mathcal{W}_1\cap\mathcal{V}(\Sigma) \neq \{0\}$.
\end{proof}

\subsection{Proof of Theorem~\ref{thm: windowedDetectionOptimal}}\label{sect: windowedDetectionOptimal}
\begin{proof}
\textit{(If)} Let $Y(T) = \mathcal{O}_T x(0)$ and $y_\Omega = \Omega x(0)$ for some $x(0) \in \mathbb{R}^n$. Then, by construction of $\widehat{Y}(k)$,
\begin{equation}\label{eqn: optimalProof1}
	\widehat{Y}(k) = \mathcal{K}(k) A^{k-l+1}x(0).
\end{equation}
for all $k = l-1, l, \dots, T$, which means that
\begin{equation}\label{eqn: optimalProofProjectionComp}
	\Pi_{\mathcal{K}(k)}\widehat{Y}(k) = \widehat{Y}(k),
\end{equation}
for all $k = l-1, l, \dots, T$. Thus, 
\[\psi\left(\widehat{Y}(l-1)\right) = \psi\left(\widehat{Y}(l)\right) = \cdots = \psi\left(\widehat{Y}(T)\right) = \text{No Attack}.\]

\textit{(Only If)} We resort to induction. 

\underline{Base Case}: In the base case, we show that if 
\[{\psi}\left(\widehat{Y}(l-1)\right) = {\psi}\left(\widehat{Y}(l)\right) = \text{No Attack},\] 
then $Y(l) = \mathcal{O}_l x(0)$ and $y_\Omega = \Omega x(0)$ for some $x(0) \in \mathbb{R}^n$. Since ${\psi}\left(\widehat{Y}(l-1)\right) = \text{No Attack}$, we have 
\begin{equation}\label{eqn: optimalProof2}
	\widehat{Y}(l-1) = \Pi_{\mathcal{K}(l-1)}\widehat{Y}(l-1),
\end{equation} which means that
\begin{align}
	\widehat{Y}(l-1) &= \mathcal{K}(l-1) x(0),\label{eqn: optimalProof2a}\\
			&= \left[\begin{array}{c} \Omega \\ \mathcal{O}_{l-1}\end{array}\right] x(0),\label{eqn: optimalProof2b}
\end{align}
for some $x(0) \in \mathbb{R}^n$. Since $\psi\left(\widehat{Y}(l)\right) = \text{No Attack}$, we have
\begin{equation}\label{eqn: optimalProof3}
	\widehat{Y}(l) = \mathcal{O}_{l-1} x'(0).
\end{equation}
for some $x'(0) \in \mathbb{R}^n$. From equation~\eqref{eqn: optimalProof2b}, we have
\begin{equation}\label{eqn: optimalProof4}
	\left[\begin{array}{ccc} y(1)^T & \cdots & y(l-1)^T\end{array}\right]^T = \mathcal{O}_{l-2} Ax(0),
\end{equation}
and from equation~\eqref{eqn: optimalProof3}, we have
\begin{equation}\label{eqn: optimalProof5}
	\left[\begin{array}{ccc} y(1)^T & \cdots & y(l-1)^T\end{array}\right]^T = \mathcal{O}_{l-2} x'(0).
\end{equation}
The pair $(A, C)$ is observable and $l \geq n+1$, so the matrix $\mathcal{O}_{l-2}$ is injective. Thus, combining equations~\eqref{eqn: optimalProof4} and ~\eqref{eqn: optimalProof5}, we have $x'(0) = Ax(0)$. By definition of $\widehat{Y}(l)$ and substituting $x'(0) = Ax(0)$ into equation~\eqref{eqn: optimalProof3}, we have that $y(l) = CA^l x(0)$. Note that $Y(l) = \left[\begin{array}{cc} \overline{Y}(l-1)^T & y(l)^T \end{array}\right]^T$. Thus, $Y(l) = \mathcal{O}_l x(0)$ and $y_\Omega = \Omega x(0)$ for some $x(0) \in \mathbb{R}^n$. 

\underline{Induction Step}: In the induction step, we assume that if 
\[\psi\left(\widehat{Y}(l-1)\right) = \cdots = \psi\left(\widehat{Y}(T-1)\right) = \text{No Attack},\] 
then $Y(T-1) = \mathcal{O}_{T-1} x(0)$ and $y_\Omega = \Omega x(0)$ for some $x(0)~\in~\mathbb{R}^n$. We show that if $\psi\left(\widehat{Y}(T)\right) = \text{No Attack}$ as well, then $Y(T) = \mathcal{O}_T x(0)$ and $y_\Omega = \Omega x(0)$ for some $x(0) \in \mathbb{R}^n$.

Since $\psi\left(\widehat{Y}(T)\right) = \text{No Attack}$, we have
\begin{equation}\label{eqn: optimalProofInduction2}
	\widehat{Y}(T) = \mathcal{O}_{l-1} x'(0),
\end{equation}
for some $x'(0)\in \mathbb{R}^n$. From the induction hypothesis, we have that $Y(T-1) = \mathcal{O}_{T-1} x(0)$, which means that
\begin{equation}\label{eqn: optimalProofInduction3}
	\left[\!\!\!\begin{array}{ccc} y(T-l+1)^T \!\!\!& \cdots \!\!\!& Y(T-1)^T \end{array}\!\!\!\right]^T\!\!\!\!\! = \!\mathcal{O}_{l-2} A^{T-l+1} x(0).
\end{equation}
From equation~\eqref{eqn: optimalProofInduction2}, we have
\begin{equation}\label{eqn: optimalProofInduction4}
	\left[\begin{array}{ccc} y(T-l+1)^T & \cdots & Y(T-1)^T \end{array}\right]^T = \mathcal{O}_{l-2} x'(0). 
\end{equation}
The pair $(A, C)$ is observable and $l \geq n+1$, so the matrix $\mathcal{O}_{l-2}$ is injective. As a result, $x'(0) = A^{T-l+1} x(0)$. Substituting $\theta ' = A^{T-l+1}$ into equation~\eqref{eqn: optimalProofInduction4}, we have $y(T) = CA^T x(0)$. Note that $Y(T) = \left[\begin{array}{cc} Y(T-1)^T & y(T)^T \end{array}\right]^T$. Thus, $Y(T) = \mathcal{O}_T x(0)$ and $y_\Omega = \Omega x(0)$ for some $x(0) \in~\mathbb{R}^n$.
\end{proof}

\section{Numerical Example}\label{sect: numExamples}
We illustrate our results with an example of a remotely piloted aircraft subject to both nonzero state inducing attacks and zero state inducing attacks. Reference~\cite{Linehan} provides a numerical model of the longitudinal dynamics of a remotely piloted aircraft that accounts for the aircraft's physical parameters.  We describe the longitudinal dynamics of the aircraft using four state variables: horizontal velocity ($x_1$), vertical velocity ($x_2$), pitch rate ($x_3$), and pitch angle ($x_4$). The aircraft we consider has two actuators: the elevator ($u_1$) and the thrust ($u_2$). The aircraft also has three sensors: the horizontal velocity sensor ($y_1$), the vertical velocity sensor ($y_2$), and the pitch angle sensor ($y_3$).

The evolution of the state variables $x_1, \dots, x_4$ is determined by physical principles governing the longitudinal flight of the aircraft and depends on physical parameters of the aircraft such as its mass and its pitch moment. The model is linearized about an equilibrium point, so the state variables $x_1, \dots, x_4$ represent values of the internal states relative to a fixed point (e.g., $x_1$ in the linearized model is the horizontal velocity of the aircraft relative to an equilibrium horizontal velocity). The linearized, discretized model for the aircraft gives the following dynamics and sensing matrices~\cite{Linehan}:
\begin{equation}\label{eqn: exampleAMatrix}
	A = \left[\begin{array}{rrrr} 0.992 & 0.030 & -0.003 & -0.977 \\
0.025 & 0.684 & 1.847 & -0.041 \\
0.054 & -0.100 & 0.381 & -0.025 \\
0.003 & -0.006 & 0.068 & 0.999 \end{array}\right],	
\end{equation}
\begin{equation}\label{eqn: exampleCMatrix}
	C = \left[\begin{array}{cccc} 1 & 0 & 0 & 0 \\
0 & 1 & 0 & 0 \\
0 & 0 & 0 & 1 \end{array}\right].
\end{equation}
The pair $(A, C)$ in this example is observable. 

We consider an attacker modeled by the following $B$ and $D$ matrices:
\begin{equation}\label{eqn: exampleBMatrix}
	B = \left[\begin{array}{rrrr} 0.001 & 0.025 & 0& 0\\
-3.224 & -0.035 & 0& 0\\
-1.995 & -0.021 & 0 & 0 \\
-0.115 & -0.001 & 0 & 0 \end{array}\right],
\end{equation}
\begin{equation}\label{eqn: exampleDMatrix}
	D = \left[\begin{array}{cccc}0 & 0 & 1 & 0 \\
0 & 0 & 0 & 1 \\
0 & 0 & 0 & 0 \end{array}\right]. 
\end{equation}
The attacker can attack both actuators (elevator, $u_1$, and thrust, $u_2$) and the horizontal velocity ($y_1$) and vertical velocity ($y_2$) sensors. There exists a zero dynamics attack against the system $\Sigma = (A, B, C, D)$. 

In this numerical example, we compare the performance of a detector that does not use side information (i.e., the detector's side information matrix is $\Omega = 0$) and the performance of a detector that uses side information matrix
\[\Omega = \left[\begin{array}{cccc} 1 & 0 & 0 & 0\end{array}\right].\]
The detector with nontrivial side information knows the initial horizontal velocity $x_1(0)$. Both detectors are implementations of the windowed detector presented in Section~\ref{sect: detectorSideInfo}; the only difference between the use of side initial state information.

\begin{figure}[h!]
	\centering
	\includegraphics[keepaspectratio = true, scale = .6]{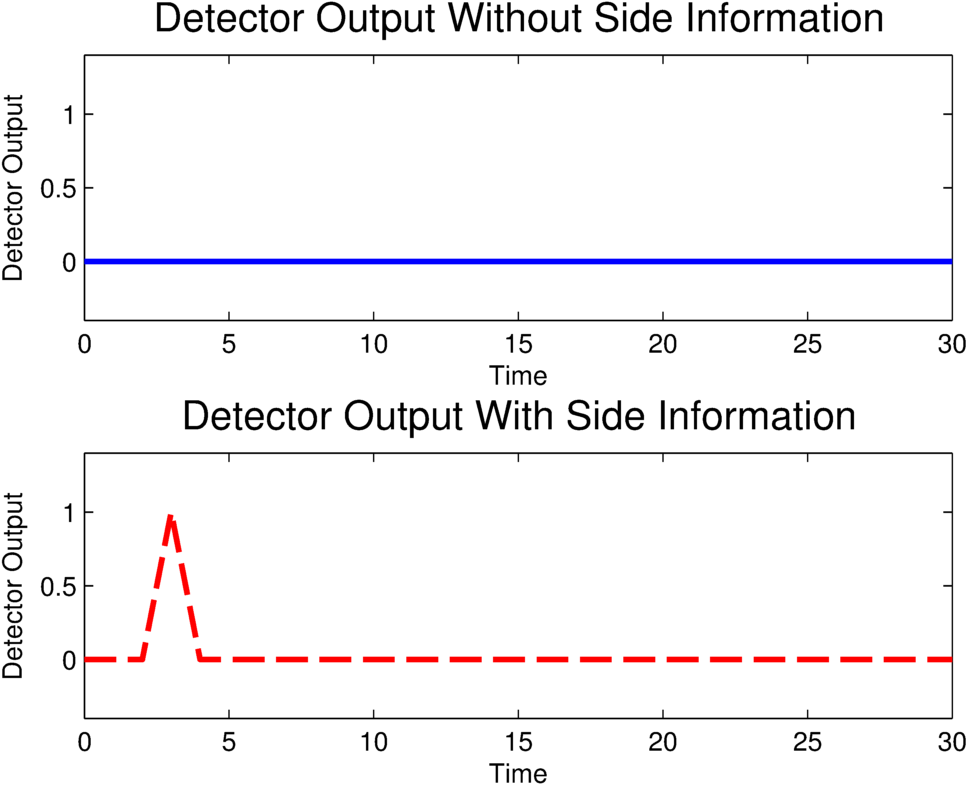}
	\caption{Detector performance without side information (top) and with side information (bottom) against zero dynamics attack.}
	\label{fig: detectorPerformance}
\end{figure}

We construct a zero dynamics attack (as defined in~\cite{TeixeiraModels} and~\cite{Pasqualetti}) against the remotely piloted aircraft. Following equation~\eqref{eqn: zeroDynamicsIndividualAttack}, we construct the zero dynamics attack component wise as
\begin{equation}\label{eqn: exampleZeroDynamicsAttack}
	a(k) \!=\! (10)(.9779)^k\left[\begin{array}{cccc} .0324 & 0 & -.6396 & .3007\end{array}\right]^T,
\end{equation}
where $k=0, \dots, 30$. The performance of the two detectors are shown in Figure~\ref{fig: detectorPerformance}.
The detector without side information is unable to detect the zero dynamics attack -- the detector outputs $0$, equivalent to ``No Attack'' for all times. The detector with side information is able to detect the zero dynamics attack -- the detector has an output of $1$, equivalent to ``Attack'' at time t= 3.

\section{Conclusion}\label{sect: conclusion}
In this paper, we studied the effect of side initial state information on the dynamic detection of data deception attacks against cyber-physical systems. First, an undetectable attack induces a state in the intersection of the system's weakly unobservable subspace, $\mathcal{V}(\Sigma)$, and the null space of the side information matrix, $\mathcal{N}\left(\Omega\right)$. Second, an undetectable attack $E(T)$ has an undetectable extension to any $T' > T$ if and only if the sum of the change in state produced by the attack, $\mathcal{C}_T E(T)$, and the zero-input state response of the state induced by the attack, $A^{T+1}\theta$, belongs to the system's weakly unobservable subspace, $\mathcal{V}(\Sigma)$. Third, there exists an arbitrarily long zero state inducing attack if and only if the intersection of the system's weakly unobservable subspace, $\mathcal{V}(\Sigma)$, and the system's output-nulling reachable subspace over one step, $\mathcal{W}_1$, is nonzero. Finally, we designed an attack detector that uses side information and detects all attacks that are not undetectable. 

\bibliography{References}
\end{document}